\documentclass[10pt]{elsarticle}
\usepackage[cp1251]{inputenc}
\usepackage{amsmath}
\usepackage{amssymb}
\usepackage{amsfonts}
\usepackage{mathtools}
\usepackage{dsfont}
\usepackage{graphicx}
\usepackage{floatflt,epsfig}
\usepackage{lineno,hyperref}
\usepackage{enumerate}
\usepackage{colortbl}
\usepackage{array,tabularx,tabulary,booktabs}
\usepackage{longtable}
\usepackage{multirow}
\usepackage{wrapfig}
\usepackage{subcaption}
\usepackage{caption}
\usepackage[table]{xcolor}
\usepackage{adjustbox}
\usepackage{rotating}
\usepackage{arydshln}
\usepackage{bm}
\usepackage{tabularray}
\newcolumntype{^}{>{\currentrowstyle}}

\usepackage{enumitem}
\usepackage{tabularx}
\usepackage{graphicx}

\usepackage{float}
\journal{Siberian Electronic Mathematical Reports}
\newtheorem{lemma}{Lemma}

\newtheorem{theorem}{Theorem}
\newtheorem{corollary}{Corollary}

\newcommand{\proof}{\medskip\noindent{\bf Proof.~}}
\usepackage{array}
\newcolumntype{P}[1]{>{\raggedright\arraybackslash}p{#1}}

\begin{document}
\renewcommand{\abstractname}{Abstract}
\renewcommand{\refname}{References}
\renewcommand{\arraystretch}{0.9}
\thispagestyle{empty}
\sloppy

\begin{frontmatter}
\title{Constructing segments of quadratic length in $Spec(T_n)$ through segments of linear length}

\author[01,02]{Artem Kravchuk}
\ead{artemkravchuk13@gmail.com}

\address[01]{Sobolev Institute of Mathematics, Ak. Koptyug av. 4, Novosibirsk 630090, Russia}
\address[02]{Novosibirsk State University, Pirogova str. 2, Novosibirsk, 630090, Russia}

\begin{abstract}
A \emph{Transposition graph} $T_n$ is defined as a Cayley graph over the symmetric group $Sym_n$ generated by all transpositions. It is known that the spectrum of $T_n$ consists of integers, but it is not known exactly how these numbers are distributed. In this paper we prove that integers from the segment $[-n, n]$ lie in the spectrum of $T_n$ for any $n\geqslant 31$. Using this fact we also prove the main result of this paper that a segment of quadratic length with respect to $n$ lies in the spectrum of $T_n$.
\end{abstract}

\begin{keyword}
Transposition graph; integral graph; spectrum; 
\vspace{\baselineskip}
\MSC[2010] 05C25\sep 05E10\sep 05E15
\end{keyword}
\end{frontmatter}

\section{Introduction}\label{sec1}

A \emph{Transposition graph} $T_n$ is defined as a Cayley graph over the symmetric group $Sym_n$ generated by all transpositions. The graph $T_n, n\geqslant 2$, is a connected bipartite  $C_n^2$-regular graph~\cite{K08}. The \emph{spectrum} of a graph is defined as a multiset of distinct eigenvalues of its adjacency matrix together with their multiplicities. It was shown in ~\cite[Lemma 3]{KY97} that $T_n$ is integral which means that its spectrum $Spec(T_n)$ consists of integers. Later and independently in ~\cite[Theorem 2]{KL20}, it was also shown that  $T_n$ is integral.  Since the Transposition graph is bipartite, its spectrum is symmetric about zero \cite{BH12}, i.e., if number $k$ lies in $Spec(T_n)$  then number $-k$ lies in $Spec(T_n)$ too. The largest eigenvalue of $T_n$ is equal to $C_{n}^2$ which implies that all eigenvalues of $T_n$ lie in the interval $[-C_n^2, C_n^2]$.

 A nonincreasing sequence $p = (n_1, \dots, n_k) \vdash n $, $k\geqslant 1$, for which $n=\sum\limits_{j=1}^kn_j$, is called a \emph{partition} of $n$.  It is follows from ~\cite[eq. 3]{KY97} that any partition $p\vdash n$ corresponds to an eigenvalue $\lambda(p) \in Spec(T_n)$ by the following expression: 

\begin{equation}\label{transp_eigen}
    \lambda(p) = \sum_{j=1}^k \frac{n_j(n_j-2j+1)}{2}. 
\end{equation}

The \emph{conjugate of partition $p$}, denoted by $p'$, is also a partition of $n$ where its parts are the nonincreasing sequence $p' = (n_1', n_2',\ldots,n_{\kappa}')$, where $n_j' = \sum\limits_{i|n_i \geqslant j}1$. This paper uses that $n_1'$ is equal to the length of $p$. It is follows from ~\cite[Lemma 8]{DS81} that 

\begin{equation}\label{conj_eigen}
    \lambda(p) = - \lambda(p'). 
\end{equation}

Note that different partitions can correspond to the same eigenvalue. For example, the partitions $(4, 1, 1)$ and $(3, 3)$ correspond to the eigenvalue $3$ of $T_6$.

Thus, if one enumerates all partitions of number $n$ and substitutes them into expression (\ref{transp_eigen}) one can obtain all eigenvalues of $T_n$. We write $f(n) \sim g(n)$ if $\frac{f(n)}{g(n)} \rightarrow 1$ as $n \rightarrow \infty$. The number of partitions $p(n)$ has the following asymptotics: $p(n)\sim \frac{exp(\pi \sqrt{\frac{2}{3}}\sqrt{n-\frac{1}{24}})}{4n\sqrt{3}}$ \cite{HR1918}, so it is no longer possible to enumerate all partitions for large $n$. 

Moreover, looking at expression (\ref{transp_eigen}), it is even not possible to answer simple questions about $Spec(T_n)$: 
\begin{itemize}
    \item is number $0$ an eigenvalue of $T_n$?
    \item  are all integers from the interval $[0, k], k\in \mathbb{N}$, lie in the spectrum of $T_n$?
    \item what is the asymptotics for the number of unique values in the spectrum of $T_n$?
\end{itemize}
Ideally, one would like to be able to answer the question whether a given integer number $k \in [-C_n^2, C_n^2]$ lies in the spectrum of $T_n$. The following result is known.

\begin{theorem} \label{KK22} {\rm ~\cite[Theorem~2]{KK22}}
 For any integer $k\geqslant 0$, there exists $n(k)$ such that for any $n\geqslant n(k)$ and any $m \in \{0, \dots, k\}$, $m \in Spec(T_n)$. 
 \end{theorem}

In the proof of this theorem, $n(k)$ is chosen as $n(k)=10k+4$. Consequently, all integers within the range $[-\frac{n-4}{10}, \frac{n-4}{10}]$ are encompassed within the spectrum of $T_n$. The following result extends the result of Theorem \ref{KK22}.

\begin{theorem}\label{segment_0_n2_theorem}{\rm ~\cite[Theorem~3]{KK23}}
For  any $n \geqslant 19$, all integers from the segment $[-\frac{n-4}{2}, \frac{n-4}{2}]$ lie in the spectrum of $T_n$.
\end{theorem}
The following result of this paper improves Theorem \ref{segment_0_n2_theorem}.

\begin{theorem}\label{th-3}
For any $n \geqslant 31$, all integers from the segment $[-n, n]$ lie in the spectrum of $T_n$.
\end{theorem}

Theorems \ref{KK22}, \ref{segment_0_n2_theorem}, \ref{th-3} have in common that all integers in some neighborhood $[-x, x]$ of zero lie in the spectrum of $T_n$. In the proof of these theorems the segment $[0, x]$ is splitted into subsegments $S=\bigcup s_i, s_i=[\underline{x_i}, \overline{x_i}]$ and individual integer numbers $A=\{a_1, \dots, a_k\}$ which may depend on the parity of $n$. Moreover, $S$ and $A$ are chosen so that $S \bigcup A$ contains all integers from the segment $[0, x]$. For each segment $s_i$, there is a family of partitions such that all integers from the segment $s_i$ are covered by the eigenvalues corresponding to the partitions from this family. Similar, for each individual integer $a_i$, there is a partition $p_i$ such that $\lambda(p_i)=a_i$.

For example \cite[Lemma~4]{KK22}, for any odd $n \geqslant 7$, the partition $\bigl(\frac{n-2\lambda+1}{2}, \lambda+2, 2\times (\lambda-1), 1\times \frac{n-4\lambda-1}{2}\bigr)$ corresponds to the eigenvalue $\lambda \in \mathbb{N}$, where $\lambda \in [1,  \frac{n-3}{4}]$. The notation $(n_1, ..., n_j, 2\times t_2, 1 \times t_1)$ means that the number $2$ is repeated in the partition $t_2$ times and the number $1$ is repeated $t_1$ times.

Thus, in the proof of Theorem \ref{th-3}, we obtain a set of partitions  $P=\bigcup p_i$ such that for each $k \in [0, n], k \in \bigcup \lambda (p_i)$. 
Theorem \ref{th-3} gives the following corollary, which is useful for the proof of Theorem \ref{th-4}.

\begin{corollary} \label{cor-1}
The first part of all partitions used to prove Theorem \ref{th-3} does not exceed $\frac{n+3}{2}$.
\end{corollary}

Note that Theorem \ref{th-3} shows that all integers from the interval $[-n, n]$ lie in the spectrum of $T_n$ and the length of this interval has asymptotics $O(n)$.

The main result of this paper shows that all integers from the set $[-y_2, -y_1] \bigcup [y_1, y_2]$ lie in the spectrum of $T_n$. Moreover, the length of the segment $[y_1, y_2]$ has $O(n^2)$ asymptotics. 

\begin{theorem}\label{th-4}
    For any $n \geqslant 48$, all integers from the segments $[-y_2, -y_1]$ and $[y_1, y_2]$, lie in the spectrum of $T_n$, where $y_1=C_{\lceil \frac{n}{3} \rceil + 1}^2 - 2(\lfloor \frac{2n}{3}\rfloor - 1)$ and $y_2=C_{\lfloor \frac{2n + 1}{3} \rfloor}^2$.
\end{theorem}
The paper is structured as follows. First, in Section \ref{sec3} we prove Theorem \ref{th-3}. The proof of Theorem~\ref{th-3} is based on several technical lemmas, which are proved in Section~\ref{sec-tech}. After proving Theorem~\ref{th-3}, we use it and Corollary~\ref{cor-1} to prove Theorem~\ref{th-4} in Section~\ref{sec2}. 

\section{Proof of Theorem~\ref{th-3}}\label{sec3}
In what follows further the notation $S\in Spec(T_n)$ where $S$ is a segment $[\underline{x}, \overline{x}]$ means that all integers from this segment are eigenvalues of $T_n$. Similar, for finite set of numbers $A$ we say that  $A\in Spec(T_n)$ if all elements of this set lie in the spectrum of $T_n$.

To prove the theorem, we  use the approach described in \cite{KK23}. We split the segment $[0, n]$ into subsegments $S_1, S_2$ and into two sets $A_1$, $A_2$ of numbers  depending on the parity of $n$ such that $S_1\cup S_2 \cup A_1 \cup A_2$ contains all integers of the interval $[0, n]$. The segments $S_1, S_2$ and the sets $A_1, A_2$ are shown in Table \ref{split-0n}.

\renewcommand{\arraystretch}{2}

\begin{table}[!h]
\centering
\begin{tabularx}{0.95\textwidth}{|c|c|c|c|X|}
\hline
	$n$ & $S_1$  & $A_1$ & $S_2$ & $A_2$ \\
 \hline
 odd & $[0, \frac{n-1}{2}]$ & $\{\frac{n+1}{2}, \frac{n+3}{2}, \frac{n+5}{2}\}$ & $[\frac{n+7}{2}, n-7]$& \multirow{2}{12em}{$\{n-6, n-5, \dots, n\}$} \\

 even & $[0, \frac{n-4}{2}]$ & $\{\frac{n-2}{2}, \frac{n}{2}, \frac{n+2}{2}\}$ & $[\frac{n+4}{2}, n-6]$ & \\
 \hline
\end{tabularx}
    \caption{ \label{split-0n}  Splitting of segment $[0, n]$.}
\end{table}

Next, we show that $S_1 \in Spec(T_n), S_2 \in Spec(T_n)$, $A_1 \in Spec(T_n)$ and $A_2 \in Spec(T_n)$. To show this, we prove the following technical lemmas. 

\begin{lemma} \label{s_1_lemma}
    $S_1 \in Spec(T_n)$ for any $n \geqslant 19$.
\end{lemma}
\begin{corollary}\label{cor-lemma-s1}
    The first part of all partitions together with their conjugates used to cover $S_1$ does not exceed  $\frac{n + 1}{2}$.
\end{corollary}
\begin{lemma}\label{lemma-s2}
    $S_2 \in Spec(T_n)$ for any $n \geqslant 20$.
\end{lemma}

\begin{corollary}\label{cor-lemma-s2}
The first part in the partitions together with their conjugates  used in the proof of Lemma \ref{s_1_lemma} does not exceed $\frac{n + 2}{2}$.    
\end{corollary}

\begin{lemma}\label{lemma-a1}
    $A_1 \in Spec(T_n)$ for any $n \geqslant 31$. 
\end{lemma}
\begin{corollary}\label{cor-lemma-a1}
    The first part in the partitions together with their conjugates  used in the proof of Lemma \ref{lemma-a1} does not exceed $\frac{n+2}{2}$.
\end{corollary}
\begin{lemma}\label{lemma-a2}
    $A_2 \in Spec(T_n)$, for any $n\geqslant19$.
\end{lemma}
\begin{corollary}\label{cor-lemma-a2}
    The first part in the partitions together with their conjugates used in the proof of Lemma \ref{lemma-a2} does not exceed $\frac{n+3}{2}$.
\end{corollary}

The proof of these technical lemmas can be found in Section~\ref{sec-tech}. After proving the lemmas, it is directly deduced that the segment $[0, n] \in Spec(T_n)$ for any $n\geqslant 31$. Since $T_n$ is bipartite, we have $[-n, n] \in Spec(T_n)$. Since Corollaries \ref{cor-lemma-s1}-\ref{cor-lemma-a2} consider partitions together with their conjugates, it follows from (\ref{conj_eigen}) that the first part of the partitions covered by the segment $[-n, n]$ does not exceed $\frac{n+3}{2}$, which proves Corollary \ref{cor-1}.

\section{Proof of Theorem~\ref{th-4}}\label{sec2}



 Let $p=(n_1, n_2, \dots, n_k) = (n_1, p_1)$, where $p_1=(n_2, \dots, n_k)$. It turns out that the eigenvalue $\lambda(p)$ corresponding to the partition $p$ can be compactly expressed through the eigenvalue $\lambda(p_1)$ corresponding to the partition $p_1$ and the number $n_1$ using the following lemma.

\begin{lemma}
\begin{equation}\label{eigen_expr}
\lambda(p)=C_{n_1}^2 + \lambda(p_1) - (n-n_1).  
\end{equation}
   
\end{lemma}
\proof
    We prove the lemma by a direct substitution of partition $p=(n_1, n_2, \dots, n_k)$ into the expression (\ref{transp_eigen}). 
    $$\lambda(p) = \frac{1}{2}\sum_{j=1}^kn_j(n_j-2j+1) = \frac{n_1(n_1 - 1)}{2} + \frac{1}{2}\sum_{j=2}^{k}n_j(n_j-2j+1) = $$
    $$=C_{n_1}^2 + \frac{1}{2}\sum_{j=2}^{k}n_j(n_j-2(j-1) + 1 - 2) =$$ 
    $$= C_{n_1}^2 + \sum_{j=2}^{k}n_j(n_j-2(j-1)+1) - \sum_{j=2}^{k}n_j = $$
    $$=C_{n_1}^2 + \lambda(p_1) - (n-n_1).$$
\hfill $\square$
\newline

Let us consider the expression~(\ref{eigen_expr}). By varying $p_1$, one can obtain all possible eigenvalues for the partitions that have $n_1$ in the first part. 
The following conditions are imposed on $p_1$:
\begin{enumerate}[label=\arabic*.]
    \item  $p_1 \vdash n-n_1$; \label{enum1}
    \item the first part of $p_1$ must not exceed $n_1$. \label{enum2}
\end{enumerate}

Let us denote the set of partitions that satisfy these two conditions by $\mathbb{P}_{n_1}$  and the set of eigenvalues that correspond to these partitions by $\bm{\lambda_{n_1}}$. If it is proved that all integers from the interval $[-(n-n_1), n-n_1]$ lie in the $\bm{\lambda_{n_1}}$, then it follows directly from (\ref{eigen_expr}) that $[C_{n_1}^2 - 2(n-n_1), C_{n_1}^2] \in Spec(T_n)$. Note that when $n_1 \leqslant \frac{2n+1}{3}$ then $2(n-n_1)\geqslant n_1-1= C_{n_1}^2 - C_{n_1 - 1}^2$. Thus, when $n_1 \leqslant \frac{2n+1}{3}$, then $[C_{n_1}^2 - 2(n-n_1), C_{n_1}^2] \in Spec(T_n)$ implies  $[C_{n_1-1}^2, C_{n_1}^2] \in Spec(T_n)$.

Consider the segment $l=[k_1, k_2]$, where $k_1\in \mathbb{N}$ and $k_2=\lfloor\frac{2n+1}{3}\rfloor$. If for each natural $n_1 \in l$ it holds that all integers from the interval $[-(n-n_1), n-n_1]$ lie in the $\bm{\lambda_{n_1}}$, then it follows that  $[C_{k_1}^2 - 2(n-k_1), C_{k_2}^2] \in Spec(T_n)$. 

By Theorem \ref{th-3}, we have that if $n-n_1 \geqslant 31$, then $[-(n-n_1), n-n_1] \in Spec(T_{n-n_1})$. Moreover, by Corollary \ref{cor-1}, to cover the segment $[-(n-n_1), n-n_1]$ we are able to choose such partitions that their first part does not exceed $\frac{n-n_1 + 3}{2}$. Inequality $n_1 \geqslant \frac{n-n_1 + 3}{2}$ can be rewritten as $n_1 \geqslant \frac{n}{3} + 1$. So, if $n_1 \geqslant \frac{n}{3} + 1$ and $n-n_1 \geqslant 31$, then $[-(n-n_1), (n-n_1)] \in \bm{\lambda_{n_1}}$.

Thus, for any natural $n_1 \in [\frac{n}{3} + 1, \frac{2n+1}{3}]$ and $n-n_1 \geqslant 31$ we have that $[C_{n_1-1}^2 - 2(n-n_1), C_{n_1}^2]\in Spec(T_n)$. Therefore, $[C_{\lceil \frac{n}{3} \rceil + 1}^2 - 2(n - \lceil \frac{n}{3} \rceil - 1), C_{\lfloor \frac{2n + 1}{3} \rfloor}^2] \in Spec(T_n)$ when $n-n_1 \geqslant 31$. Note that $n - \lceil\frac{n}{3}\rceil - 1=\lfloor \frac{2n}{3}\rfloor - 1$, hence $[C_{\lceil \frac{n}{3} \rceil + 1}^2 - 2(\lfloor \frac{2n}{3}\rfloor - 1), C_{\lfloor \frac{2n + 1}{3} \rfloor}^2] \in Spec(T_n)$. Moreover, $[-C_{\lfloor \frac{2n + 1}{3} \rfloor}^2, -([C_{\lceil \frac{n}{3} \rceil + 1}^2 - 2(\lfloor \frac{2n}{3}\rfloor - 1))] \in Spec(T_n)$ too, since $T_n$ is bipartite.
The inequalities $n - n_1 \geqslant 31$ and $n_1\geqslant \frac{n}{3} + 1$ imply that $n \geqslant 48$, which completes the proof.

\hfill $\square$


\section{Proof of technical lemmas}\label{sec-tech}
Let the partition $p$ is given by $p=(p_l, 2\times t_2, 1\times t_1)$, where $p_l$ is a sequence of length $l$. For example, a partition $(5, 4, 4, 2, 2, 2, 1, 1)$ can be represented as $(p_3, 2\times 3, 1\times 2)$, where $p_3=(5, 4, 4)$. Note that for a single partition the representation in this form is not the only one. For example, we can represent the partition $(5, 4, 4, 2, 2, 2, 1, 1)$ as $(p_4, 2\times 2, 1\times2)$, where $p_4=(5, 4, 4, 2)$. Therefore, in what follows we emphasize how $p_l$ and $t_1, t_2$ are chosen.
\begin{lemma}\label{lemma-2-1}
If $p=(p_l, 2\times t_2, 1\times t_1)$, then
$$\lambda(p) = \lambda(p_l)-f(t_1, t_2, l) = \lambda(p_l)-  \underbrace{(t_2-1)^2 - \frac{1}{2}(t_1-\frac{1}{2})^2 - t_1t_2 - l(2t_2 + t_1)+\frac{9}{8}}_{f(t_1, t_2, l)}.$$
\end{lemma}

\proof
By (\ref{transp_eigen}) we have:

$$
\lambda(p) = \lambda(p_l) + \frac{1}{2}\sum\limits_{j=l+1}^{l+t_2} 2(3-2j) + \frac{1}{2}\sum\limits_{l+t_2+1}^{l+t_1+t_2}(2-2j)=
$$
$$
=\lambda(p_l) + 3t_2 + t_1 - 2\sum\limits_{j=l+1}^{l+t_2}j - \sum\limits_{l+t_2+1}^{l+t_1+t_2}j=
$$

$$
=\lambda(p_l) + 3t_2 + t_1 - (2l+t_2+1)t_2 - \frac{2l + 2t_2 + t_1 + 1}{2}t_1=
$$

$$
= \lambda(p_l) - (t_2-1)^2 - \frac{1}{2}(t_1-\frac{1}{2})^2 - t_1t_2 - l(2t_2 + t_1)+\frac{9}{8}.
$$

\hfill $\square$

\begin{corollary}\label{cor-2}
    If $p=(p_l, 1\times t_1)$, then
    \begin{equation}\label{t1-expr}
        \lambda(p) = \lambda(p_l) - f(l, t_1) = \lambda(p_l) - \underbrace{t_1(\frac{t_1 - 1}{2} + l)}_{f(t_1, l)}.
    \end{equation}
\end{corollary}

\medskip\noindent{\bf Proof of Lemma~\ref{s_1_lemma}.}  The proof is straightforward from the technical lemmas of \cite{KK23}. For the reader's convenience, we have presented the summarization of these lemmas in Table~\ref{kk23-table}.

\medskip\noindent{\bf Proof of Corollary~\ref{cor-lemma-s1}.} 

The proof if straightforward from Table \ref{kk23-table}. The largest value of the first part is $\frac{n+1}{2}$ and is obtained on the partition $(\frac{n+1}{2}, 1\times \frac{n-1}{2})$. The largest value among the conjugate partitions is $\frac{n+1}{2}$ and is achieved on the partition conjugate to $(\frac{n+1}{2}, 1\times \frac{n-1}{2})$. Note that the partition $(\frac{n+1}{2}, 1\times \frac{n-1}{2})$ is self-conjugate and it follows from (\ref{conj_eigen}) that if the partition $p$ is self-conjugate, then $\lambda(p) = 0$.

\medskip\noindent{\bf Proof of Lemma~\ref{lemma-s2}.} We prove this lemma by considering four cases, for each parity of $n$ and $\lambda$. For each parity of $n$, we choose two families of partitions $p_{\lambda}$ that depend on the parity of $\lambda$. Then we show that eigenvalues corresponding to the chosen families of partitions cover all integers from the segment $S_2$ which have the same parity as $\lambda$.

\medskip\noindent{\bf Case 1,} $n$ is odd, $\lambda$ is odd:

 $p_{\lambda} = \bigl(\underbrace{\frac{\lambda}{2} + \frac{3}{2}, \frac{n + 2}{2} - \frac{\lambda}{2}, 3}_{p_l, l=3}, 2\times\underbrace{(\frac{n-4}{2} - \frac{\lambda}{2})}_{t_2}, 1\times \underbrace{(\lambda - \frac{n + 3}{2})}_{t_1}\bigr)$.

By~(\ref{transp_eigen}) we have:
$$\lambda(p_l)=\frac{1}{8}n^2 - \frac{1}{4}n\lambda - \frac{1}{4}n + \frac{1}{4}\lambda^2 + \frac{3}{4}\lambda - \frac{29}{8}.$$

In addition, by ~(\ref{lemma-2-1}) it follows that:

$$f(t_1, t_2, l)=\frac{1}{8}n^2 - \frac{1}{4}n\lambda - \frac{1}{4}n + \frac{1}{4}\lambda^2 - \frac{1}{4}\lambda - \frac{29}{8}.$$

Finally, using ~(\ref{lemma-2-1}) we have $\lambda(p) = \lambda(p_l) + f(t_1, t_2, l) = \lambda$.
Partition $p$ holds for any odd integer $\lambda \in [\frac{n+3}{2}, n-4]$. Note that the interval $[\frac{n+3}{2}, n-4]$ makes sense for any $n \geqslant 11$. Also notice that segment $[\frac{n+7}{2}, n-7]$ lies inside segment $[\frac{n+3}{2}, n-4]$.

\medskip\noindent{\bf Case 2,} $n$ is odd, $\lambda$ is even:

$p_{\lambda}=\bigl(\underbrace{\frac{\lambda}{2}, \frac{n+3}{2} - \frac{\lambda}{2}, 5}_{p_l, l=3},
2\times\underbrace{(\frac{n-3}{2} - \frac{\lambda}{2})}_{t_2}, 1\times\underbrace{(\lambda -\frac{n + 7}{2})}_{t_1}\bigr)$.

For $\lambda(p_l)$ and $f(t_1, t_2, l)$ we have the following expressions:

$$\lambda(p_l) = \frac{1}{8}n^2 - \frac{1}{4}n\lambda  + \frac{1}{4}\lambda^2 - \frac{1}{4}\lambda - \frac{9}{8};$$

$$f(t_1, t_2, l)=\frac{1}{8}n^2 - \frac{1}{4}n\lambda  + \frac{1}{4}\lambda^2 - \frac{5}{4}\lambda - \frac{9}{8}.$$

Finally, $\lambda(p) = \lambda(p_l) + f(t_1, t_2, l) = \lambda$. Note that $p$ holds when even $\lambda \in [\frac{n+7}{2}, n-7]$, and $[\frac{n+7}{2}, n-7]$ holds for any $n \geqslant 21$.

\medskip\noindent{\bf Case 3,} $n$ is even, $\lambda$ is odd:

$p_{\lambda}=\bigl(\underbrace{\frac{\lambda}{2}+\frac{3}{2}, \frac{n+3}{2}-\frac{\lambda}{2}}_{p_l,l=2}, 2\times\underbrace{(\frac{n-1}{2}-\frac{\lambda}{2})}_{t_2}, 1 \times \underbrace{(\lambda - \frac{n+4}{2})}_{t_1}\bigr)$.

The following holds:

$$\lambda(p_l)=\frac{1}{8}n^2 - \frac{1}{4}n\lambda  + \frac{1}{4}\lambda^2 + \frac{1}{2}\lambda - \frac{3}{4};$$

$$f(t_1, t_2, l)=\frac{1}{8}n^2 - \frac{1}{4}n\lambda  + \frac{1}{4}\lambda^2 - \frac{1}{2}\lambda - \frac{3}{4}.$$
Finally, $\lambda(p) = \lambda(p_l) + f(t_1, t_2, l) = \lambda$. The following constraints are imposed on $\lambda$ and $n$,  so that $p$ holds: odd $\lambda \in [\frac{n+4}{2}, n-1]$, $n \geqslant 6$. Also notice that segment $[\frac{n + 4}{2}, n - 6]$ lies inside segment $[\frac{n+4}{2}, n-1]$.

 \medskip\noindent{\bf Case 4,} $n$ is even, $\lambda$ is even:

 $p_{\lambda}=\bigl(\underbrace{\frac{\lambda}{2}+1, \frac{n+2}{2}-\frac{\lambda}{2}, 4}_{p_l, l=3},  2\times\underbrace{(\frac{n-4}{2}-\frac{\lambda}{2})}_{t_2}, 1 \times \underbrace{(\lambda - \frac{n+4}{2})}_{t_1}\bigr).$

The following holds:

$$\lambda(p_l)=\frac{1}{8}n^2 - \frac{1}{4}n\lambda -\frac{1}{4}n + \frac{1}{4}\lambda^2 + \frac{1}{2}\lambda - 3;$$

$$f(t_1, t_2, l)=\frac{1}{8}n^2 - \frac{1}{4}n\lambda -\frac{1}{4}n  + \frac{1}{4}\lambda^2 - \frac{1}{2}\lambda - 3. $$

 Note that $p$ holds for even $\lambda \in [\frac{n + 4}{2}, n - 6]$ and for any $n \geqslant 16$.

 At the end of the proof, we would like to mention that for any $n\geqslant20$ the restrictions for all four cases hold.

\medskip\noindent{\bf Proof of Corollary \ref{cor-lemma-s2}.}
We need to consider a family of partitions in each of the four cases from Lemma~\ref{lemma-s2}. For each family, we have found an interval in which $\lambda$ can lie. Let us substitute the maximum value from this interval into the partitions and choose the maximum value at the first position in these four cases. The maximum value is reached in Case 3 and is equal to $\frac{n+2}{2}$.  Moreover, the maximum value of the first part for the conjugate partition is reached in Case 3 too and equal to $\frac{n-2}{2}$.

\medskip\noindent{\bf Proof of Lemma \ref{lemma-a1}.}
For each $a \in A_1$ we give a particular partition $p \vdash n$ such that $\lambda(p) = a$. All the given partitions have the form $(p_l, 1 \times l)$, so we use (\ref{t1-expr}) to prove that $\lambda(p)=a$. For example, for $\frac{n+1}{2}$ in the case of odd $n$, we take the partition $p=\bigl(\underbrace{\frac{n-11}2, 7, 4,3,3}_{p_l, l=5}, 1\times\underbrace{(\frac{n-23}2)}_{t_1}\bigr) \vdash n$. For this partition we have:

$$\lambda(p_l) = \frac{1}{2}\bigl(\frac{n-11}{2}\cdot\frac{n-13}{2} + 28 - 4 - 12 -18 \bigr) = \frac{1}{8}n^2 - 3n + \frac{119}{8},$$

$$f(t_1, l) = \frac{n-23}{2}\bigl(\frac{n-25}{4} + 5\bigr) = \frac{n^2}{8} - \frac{7}{2}n + \frac{115}{8}.$$

Then, by (\ref{t1-expr}) we get $\lambda(p) = \lambda(p_l) - f(t_1, l) = \frac{n+1}{2}$.  Note that $p$ holds for any $n\geqslant 25$.  The given $\lambda(p_l), f(t_1, l), \lambda(p)$ together with the constraints on $n$ for the eigenvalue $\frac{n+1}{2} \in A_1$ can be found in the first row of Table \ref{proof-a1}.

The partitions $p \vdash n$ along with $\lambda(p_l)$ and $f(t_1, l)$ for the remaining numbers from $A_1$ are summarized in Table \ref{proof-a1}. In addition, Table \ref{proof-a1} shows the limitations under which the partition $p$ holds. Note that for $n=31$ all limitations on $n$ listed in Table~\ref{proof-a1} are satisfied.
$\square$
\hfill

 \begin{table}[!h]
 \centering
 \small
 \setlength{\tabcolsep}{2.5pt}
\begin{tabular}{|c|c|c|c|c|}
	\hline
	  $p$ & $\lambda(p_l)$ & $f(t_1, l)$ & $\lambda(p)$ &  limitations \\
	\hline
   \begin{tabular}{c}
      $\bigl(\underbrace{\frac{n-11}2, 7, 4,3,3}_{p_l, l=5}, 1\times\underbrace{(\frac{n-23}2)}_{t_1}\bigr)$
      \\ $\bigl(\underbrace{\frac{n-1}2, 4}_{p_l, l=2}, 1\times\underbrace{(\frac{n-7}2)}_{t_1}\bigr)$ 
      \\ $\bigl(\underbrace{\frac{n-3}2, 5, 2}_{p_l, l=3}, 1\times\underbrace{(\frac{n-11}2)}_{t_1}\bigr)$
 \end{tabular} & \begin{tabular}{c}
      $\frac{1}{8}n^2 - 3n + \frac{119}{8}$  \\
      $\frac{1}{8}n^2 - \frac{1}{2}n + \frac{19}{8}$ \\
    $\frac{1}{8}n^2 - n + \frac{31}{8}$
 \end{tabular}
 
 & \begin{tabular}{c}
      $\frac{1}{8}n^2 - \frac{7}{2}n + \frac{115}{8}$  \\
      $\frac{1}{8}n^2 - n + \frac{7}{8}$ \\ 
      $\frac{1}{8}n^2 - \frac{3}{2}n + \frac{11}{8}$
 \end{tabular} 
 
 & \begin{tabular}{c} 
 $\frac{n+1}{2}$ \\ $\frac{n+3}{2}$ \\ $\frac{n+5}{2}$
 \end{tabular}
 & \begin{tabular}{c}
      $odd \ n\geqslant 25$  \\
      $odd \ n\geqslant 9$ \\
      $odd \ n \geqslant 13$
      
 \end{tabular}
\\
\hline
 \begin{tabular}{c}
      $\bigl(\underbrace{\frac{n-16}2, 8, 5, 4, 3, 3}_{p_l,l=6}, 1\times\underbrace{(\frac{n-30}2)}_{t_1}\bigr)$
      \\ $\bigl(\underbrace{\frac{n+2}2}_{p_l, l=1}, 1\times\underbrace{(\frac{n-2}2)}_{t_1})$ 
      \\ $(\underbrace{\frac{n-12}2, 8, 3, 3, 3, 2}_{p_l,l=6}, 1\times\underbrace{(\frac{n-26}2)}_{t_1}\bigr)$
 \end{tabular} & 
 \begin{tabular}{c}
      $\frac{1}{8}n^2 - \frac{17}{4}n + 29$  \\
      $\frac{1}{8}n^2 + \frac{1}{4}n$ \\ 
      $\frac{1}{8}n^2 - \frac{13}{4}n + 14$
 \end{tabular}
 & 
 \begin{tabular}{c}
      $\frac{1}{8}n^2 - \frac{19}{4}n + 30$  \\
      $\frac{1}{8}n^2 - \frac{1}{4}n$ \\
      $ \frac{1}{8}n^2 - \frac{15}{4}n + 13$
 \end{tabular}
 & \begin{tabular}{c} 
 $\frac{n-2}{2}$ \\ $\frac{n}{2}$ \\ $\frac{n+2}{2}$
 \end{tabular}  & 
 \begin{tabular}{c}
      $even \ n \geqslant 32$ \\
      $even \ n \geqslant 2$ \\
      $even \ n \geqslant 28$
 \end{tabular}
\\
\hline
	\end{tabular}	
    \caption{ \label{proof-a1}  Partitions $p\vdash n$ and calculation of $\lambda(p)$ for them for $A_1$.}
 \end{table}

\medskip\noindent{\bf Proof of Corollary \ref{cor-lemma-a1}.} The proof is straightforward from the first column of Table~\ref{proof-a1}. The maximum value of the first part among all partitions is $\frac{n+2}{2}$ and is obtained on the partition $\bigl(\frac{n+2}{2}, 1\times (\frac{n-2}{2})\bigr)$. For conjugate partitions the maximum is obtained on the partition which is conjugate to the partition $\bigl(\frac{n+2}{2}, 1\times (\frac{n-2}{2})\bigr)$ and is equal to $\frac{n}{2}$.

\medskip\noindent{\bf Proof of Lemma \ref{lemma-a2}.} We prove this lemma in a similar manner as Lemma \ref{lemma-a1}. For every $a \in A_2$ and every parity of $n$, we give a partition $p \vdash n$ such that $\lambda(p) = a$. For example, if $a=n$ and $n$ is odd we take the partition $p=\bigl(\underbrace{\frac{n+3}2}_{p_l, l=1}, 1\times\underbrace{(\frac{n-3}2)}_{t_1}\bigr) \vdash n$, for which we have:

$$\lambda(p_l) = \frac{1}{2}\cdot \frac{n+3}{2}\cdot \frac{n+1}{2} = \frac{1}{8}n^2 + \frac{1}{2}n + \frac{3}{8},$$

$$f(t_1, l) =\frac{n-3}{2}\cdot\frac{n-1}{4} = \frac{1}{8}n^2 - \frac{1}{2}n + \frac{3}{8}.$$

By (\ref{t1-expr}) we have that $\lambda(p) = \lambda(p_l) - f(t_1, l) = n$. Note that $p$ holds for any $ n\geqslant 3$. The given $\lambda(p_l), f(t_1, l), \lambda(p)$ together with the constraints on $n$ for the odd $n \in A_2$ are located in the first row of Table \ref{proof-a2}.

The partitions $p \vdash n$ along with $\lambda(p_l)$ and $f(t_1, l)$ for the remaining numbers from $A_2$ are summarized in Table \ref{proof-a2}. This table is divided into cells for each $a \in A_2$. The top row of each cell shows the partition for odd $n$, and the bottom row shows the partition for even $n$. In addition, Table \ref{proof-a2} shows the limitations under which the partition $p$ holds. Note that for $n=19$ all limitations on $n$ listed in the table are satisfied.
$\square$
\hfill

\begin{table}[!h]
\small
\begin{tabular}{|c|c|c|c|c|}
	\hline
	  $p$ & $\lambda(p_l)$ & $f(t_1, l)$ & $\lambda(p)$ & $n$ limitations \\
	\hline
 \begin{tabular}{c}
      $\bigl(\underbrace{\frac{n+3}2}_{p_l, l=1}, 1\times\underbrace{(\frac{n-3}2)}_{t_1}\bigr)$
      \\ $\bigl(\underbrace{\frac{n}2, 4}_{p_l, l=2}, 1\times\underbrace{(\frac{n}2-4)}_{t_1}\bigr)$ 
 \end{tabular} 
 & 
 \begin{tabular}{c}
      $\frac{1}{8}n^2 + \frac{1}{2}n + \frac{3}{8}$  \\
       $\frac{1}{8}n^2 - \frac{1}{4}n + 2$ 
 \end{tabular}
 & 
 \begin{tabular}{c}
      $\frac{1}{8}n^2 - \frac{1}{2}n + \frac{3}{8}$  \\
      $\frac{1}{8}n^2 - \frac{5}{4}n + 2$ 
 \end{tabular}
 & \begin{tabular}{c} 
 $n$\end{tabular}
 & \begin{tabular}{c}
      $n\geqslant 3$ \\
      $n\geqslant 8$
 \end{tabular} 
 
 \\
 \hline
 \begin{tabular}{c}
      $\bigl(\underbrace{\frac{n+1}2,3}_{p_l, l=2}, 1\times\underbrace{(\frac{n-7}2)}_{t_1}\bigr)$
      \\ $\bigl(\underbrace{\frac{n+2}2,2}_{p_l, l=2}, 1\times\underbrace{(\frac{n}2 - 3)}_{t_1}\bigr)$ 
 \end{tabular} 
 &
 \begin{tabular}{c}
      $\frac{1}{8}n^2 - \frac{1}{8}$  \\
      $\frac{1}{8}n^2 + \frac{n}{4} - 1$
 \end{tabular}
 &
 \begin{tabular}{c}
      $\frac{1}{8}n^2 - n + \frac{7}{8}$  \\
      $\frac{1}{8}n^2 - \frac{3}{4}n$
 \end{tabular}
 &
  \begin{tabular}{c} 
 $n-1$
 \end{tabular} 
  &
  \begin{tabular}{c}
       $n\geqslant 7$ \\
       $n\geqslant 6$ 
  \end{tabular}
  \\
 \hline
 \begin{tabular}{c}
      $\bigl(\underbrace{\frac{n-1}2,4, 2}_{p_l, l=3}, 1\times\underbrace{(\frac{n-11}2)}_{t_1}\bigr)$
      \\ $\bigl(\underbrace{\frac{n-8}2, 6, 5, 3}_{p_l,l=4}, 1\times\underbrace{(\frac{n}2 - 10)}_{t_1}\bigr)$ 
 \end{tabular}
 &
 \begin{tabular}{c}
      $\frac{1}{8}n^2 - \frac{n}{2} - \frac{5}{8}$  \\
      $\frac{1}{8}n^2 - \frac{9}{4}n + 13$
      
 \end{tabular}
 &
 \begin{tabular}{c}
      $\frac{1}{8}n^2 - \frac{3}{2}n + \frac{11}{8}$  \\
      $\frac{1}{8}n^2 - \frac{13}{4}n + 15$
 \end{tabular}
 &
  \begin{tabular}{c} 
 $n-2$
 \end{tabular} 
 & 
  \begin{tabular}{c}
       $n\geqslant 11$ \\
       $n\geqslant 20$ 
  \end{tabular}
  \\
 \hline
 \begin{tabular}{c}
      $\bigl(\underbrace{\frac{n+1}2, 2, 2}_{p_l,l=3}, 1\times\underbrace{(\frac{n-9}2)}_{t_1}\bigr)$
      \\ $\bigl(\underbrace{\frac{n}2, 3, 2}_{p_l, l=3}, 1\times\underbrace{(\frac{n-10}2)}_{t_1}\bigr)$ 
 \end{tabular}
 &
 \begin{tabular}{c}
      $\frac{1}{8}n^2 - \frac{33}{8}$  \\
      $\frac{1}{8}n^2 - \frac{1}{4}n - 3$
 \end{tabular}
 &
 \begin{tabular}{c}
    $\frac{1}{8}n^2 - n - \frac{9}{8}$ \\
    $\frac{1}{8}n^2 - \frac{5}{4}n$
 \end{tabular}
 &
  \begin{tabular}{c} 
 $n-3$
 \end{tabular} 
 
 & 
  \begin{tabular}{c}
       $n\geqslant 9$ \\
       $n\geqslant 10$ 
  \end{tabular}
  \\
 \hline
 \begin{tabular}{c}
      $\bigl(\underbrace{\frac{n-1}2, 3, 3}_{p_l,l=3}, 1\times\underbrace{(\frac{n-11}2)}_{t_1}\bigr)$
      \\ $\bigl(\underbrace{\frac{n-4}2, 5, 3, 2}_{p_l,l=4}, 1\times\underbrace{(\frac{n-16}2)}_{t_1}\bigr)$ 
 \end{tabular}
 &
 \begin{tabular}{c}
       $\frac{1}{8}n^2 - \frac{1}{2}n - \frac{21}{8}$ \\
      $\frac{1}{8}n^2 - \frac{5}{4}n$
 \end{tabular}
 &
 \begin{tabular}{c}
     $\frac{1}{8}n^2 - \frac{3}{2}n + \frac{11}{8}$ \\ 
     $\frac{1}{8}n^2 - \frac{9}{4}n + 4$
 \end{tabular}
 &
  \begin{tabular}{c} 
 $n-4$
 \end{tabular} 
 & 
  \begin{tabular}{c}
       $n\geqslant 11$ \\
       $n\geqslant 16$ 
  \end{tabular}
  \\
 \hline
 \begin{tabular}{c}
      $\bigl(\underbrace{\frac{n-7}2, 5, 5, 3}_{p_l,l=4}, 1\times\underbrace{(\frac{n-19}2)}_{t_1}\bigr)$
      \\ $\bigl(\underbrace{\frac{n-2}2, 4, 2, 2}_{p_l, l=4}, 1\times\underbrace{(\frac{n-14}2)}_{t_1}\bigr)$ 
 \end{tabular} 
 &
 \begin{tabular}{c}
       $\frac{1}{8}n^2 - 2n + \frac{55}{8}$ \\
       $\frac{1}{8}n^2 - \frac{3}{4}n - 5$
 \end{tabular}
 &
 \begin{tabular}{c}
       $\frac{1}{8}n^2 - 3n + \frac{95}{8}$ \\
       $\frac{1}{8}n^2 - \frac{7}{4}n$
      
 \end{tabular}
 &
  \begin{tabular}{c} 
 $n-5$
 \end{tabular} 
 & 
  \begin{tabular}{c}
       $n\geqslant 19$ \\
       $n\geqslant 14$ 
  \end{tabular}
  \\
 \hline
   
  \begin{tabular}{c}
      $\bigl(\underbrace{\frac{n-1}2, 3, 2, 2}_{p_l,l=4}, 1\times\underbrace{(\frac{n-13}2)}_{t_1}\bigr)$
      \\ $\bigl(\underbrace{\frac{n}2, 2, 2, 2}_{p_l,l=4}, 1\times\underbrace{(\frac{n-12}2)}_{t_1}\bigr)$ 
 \end{tabular}
 &
 \begin{tabular}{c}
       $\frac{1}{8}n^2 - \frac{1}{2}n - \frac{61}{8}$ \\
       $\frac{1}{8}n^2 - \frac{1}{4}n - 9$
 \end{tabular}
 &
 \begin{tabular}{c}
        $\frac{1}{8}n^2 - \frac{3}{2}n - \frac{13}{8}$ \\
        $\frac{1}{8}n^2 - \frac{5}{4}n - 3$
 \end{tabular}
 &
 \begin{tabular}{c}
 $n-6$
 \end{tabular} 
 & 
  \begin{tabular}{c}
       $n\geqslant 13$ \\
       $n\geqslant 12$ 
  \end{tabular}
  \\
 \hline
\end{tabular}	
\caption{ \label{proof-a2}  Partitions $p \vdash n$ and calculation of $\lambda(p)$ for them for $A_2$.}
\end{table}

\medskip\noindent{\bf Proof of Corollary \ref{cor-lemma-a2}.} The proof is straightforward from the first column of Table~\ref{proof-a2}. The partition with maximum first part is $(\frac{n+3}{2}, 1\times (\frac{n-3}{2}))$. The maximum of the first part among conjugate partitions is achieved by a partition conjugate to $\bigl(\frac{n+3}{2}, 1\times (\frac{n-3}{2})\bigr)$ and is equal to $\frac{n-1}{2}$.

\section{Discussions and further research}

In this paper it is shown that for any $n\geqslant 48$,  $[y_1, y_2] \in Spec(T_n)$ and $[-y_2, y_1] \in Spec(T_n)$, where $y_1=C_{\lceil \frac{n}{3} \rceil + 1}^2 - 2(\lfloor \frac{2n}{3}\rfloor - 1)$ and $y_2=C_{\lfloor \frac{2n + 1}{3} \rfloor}^2$. Note that both of these segments have length with asymptotics $O(n^2)$.  In addition, it was shown that $[-n, n] \in Spec(T_n)$ for any $n\geqslant 31$. 

Our conjecture is that the segment $[n+1, y_1]$ lies in  $Spec(T_n)$ too and we plan to study this in future works.  If this conjecture is true, it would prove the following conjecture.

\medskip\noindent{\bf Conjecture.~} There exists $n_0$ such that $[-C_{\lfloor \frac{2n + 1}{3} \rfloor}^2, C_{\lfloor \frac{2n + 1}{3} \rfloor}^2]\in Spec(T_n)$ for any $n\geqslant n_0$.

\section*{Acknowledgements}
 The author was supported by the Mathematical Center in Akademgorodok, under agreement No. 075-15-2022-281 with the Ministry of Science and High Education of the Russian Federation.

\setlength{\arrayrulewidth}{0.1mm}
\setlength{\tabcolsep}{8pt}
\renewcommand{\arraystretch}{2}

\begin{sidewaystable}[!h]
    \centering
   \begin{tabular}{ |c|c|c|c| } 
\hline
Eigenvalues, $\lambda$ & Partitions & Limitations & Proof  \\
\hline
\multirow{2}{4em}{$\ \ \ \ \ 0$} & $\bigl(\frac{n+1}{2}, 1\times\frac{n-1}{2}\bigr)$ &$n$ is odd, $n\geqslant 1$ & \multirow{2}{4em}{\rm ~\cite{KK22}, Lemma~3}\\ 
& $\bigl(\frac{n}{2}, 2, 1\times \frac{n-4}{2}\bigr)$ & $n$ is even, $n\geqslant 4$ & \\
\hline
$[1, \frac{n-3}{4}]$ & $(\frac{n-2\lambda+1}{2}, \lambda+2, 2\times (\lambda-1), 1\times \frac{n-4\lambda-1}{2})$ & $n$ is odd, $n\geqslant 7$ & \rm ~\cite{KK22}, Lemma~5\\
\hdashline
$1$ & $\bigl(\frac{n-6}{2},4,4,2,1\times \frac{n-14}{2}\bigr)$ & $n$ is even, $n \geqslant 14$ & \rm ~\cite{KK22}, Lemma~4 \\

$[2, \frac{n-4}{4}]$ & $(\frac{n-2\lambda}{2}, \lambda+2, 3, 2\times (\lambda-2), 1\times \frac{n-4\lambda-2}{2})$ & $n$ is even, $n\geqslant 12$ & {\rm ~\cite{KK23}, Lemma~4} \\

\hline
$[\frac{n+3}{4}, \frac{n-1}{2}]$ &  $\bigl(\lambda + 1, \frac{n + 3 - 2\lambda}{2}, 2\times (\frac{n-1}{2} - \lambda), 1\times (\frac{4\lambda - n - 3}{2})\bigr)$ & $n$ is odd, $n \geqslant 5$ & {\rm ~\cite{KK23}, Lemma~5}\\ 
$[\frac{n+2}{4}, \frac{n-4}{2}]$& $\bigl(\lambda + 1, \frac{n + 2 - 2\lambda}{2}, 3,  2\times (\frac{n-4}{2} - \lambda), 1\times (\frac{4\lambda - n - 2}{2})\bigr)$ & $n$ is even, $n \geqslant 10$ & {\rm ~\cite{KK23}, Lemma~6}\\ 
\hline

\multirow{4}{4em}{$ [\frac{n-3}{4}] + 1$} & 
     $\bigl(\frac{n}{4}, \frac{n}{4}, 5, 4, 2 \times (\frac{n-20}{4}), 1\bigr)$ & $n\equiv 0 \ (mod \ 4), n\geqslant 20$ & {\rm ~\cite{KK23}, Lemma~7} \\
    &  $\bigl(\frac{n + 3}{4}, \frac{n + 3}{4}, 4,  2\times (\frac{n - 13}{4}), 1\bigr)$ & $n\equiv 1 \ (mod \ 4), n \geqslant 13$ & {\rm ~\cite{KK23}, Lemma~8} \\
    & $\bigl(\frac{n + 2}{4}, \frac{n - 2}{4}, 5, 4, 2 \times (\frac{n-22}{4}), 1, 1\bigr)$ & $n\equiv 2 \ (mod \ 4), n\geqslant 22$ & {\rm ~\cite{KK23}, Lemma~9} \\
    & $\bigl(\frac{n + 1}{4}, \frac{n + 1}{4}, 5, 3, 2\times (\frac{n - 19}{4}), 1\bigr)$ & $n\equiv 3 \ (mod \ 4), n\geqslant 19$ & {\rm ~\cite{KK23}, Lemma~10} \\
    \hline
\end{tabular}  
    \caption{ \label{kk23-table}  Summary of the technical lemmas from ~\cite{KK23}.}
\end{sidewaystable}

\end{document}